\documentclass{article}
\usepackage{amsmath}[1996/11/01]
\usepackage{amssymb,amsthm,amsxtra}
\usepackage{epsfig}

\newlength{\mytopmargin}
\newlength{\myleftmargin}
\def\zz{\rlx\hbox{\small \sf Z\kern-.4em Z}}

\newcommand{\ml}{\langle}
\newcommand{\mg}{\rangle}

\begin{document}

\vspace{1cm}
\noindent
\begin{center}{   \large \bf
Pieri-type formulas for the non-symmetric Jack polynomials}
\end{center}
\vspace{5mm}

\noindent
\begin{center}
 P.J.~Forrester and D.S.~McAnally \\

\it Department of Mathematics and Statistics, \\
University of Melbourne, Parkville, Victoria
3052, Australia
\end{center}
\vspace{.5cm}
\small
\begin{quote}
In the theory of symmetric Jack polynomials the coefficients in the
expansion of the $p$th elementary symmetric function $e_p(z)$ times
a Jack polynomial expressed as a series in Jack polynomials are known
explicitly. Here analogues of this result for the non-symmetric Jack
polynomials $E_\eta(z)$ are explored. Necessary conditions for non-zero 
coefficients 
in the expansion of $e_p(z) E_\eta(z)$ as a series in non-symmetric
Jack polynomials are given.
A known expansion formula for $z_i E_\eta(z)$
is rederived by an induction procedure, and this expansion is used to
deduce the corresponding result for the expansion of
$\prod_{j=1, \, j\ne i}^N z_j \, E_\eta(z)$, and consequently
the expansion of $e_{N-1}(z) E_\eta(z)$. In the general $p$ case
the coefficients for special terms in the expansion are presented.
\end{quote}

\vspace{.5cm}
\noindent
\section{Introduction}
\setcounter{equation}{0}
Jack polynomials and Macdonald polynomials can be defined as homogeneous
multivariable orthogonal polynomials, or as eigenfunctions of a family of
commuting differential or difference operators respectively. From the latter
viewpoint these polynomials occur in the study of certain quantum many
body systems \cite{Fo94j,Ko96}. In their most basic form the polynomials are
non-symmetric, although eigenfunctions with a prescribed symmetry with
respect to interchange of coordinates are often required in application
\cite{BF97b}. The polynomials with a prescribed symmetry can be obtained
from the non-symmetric polynomials by an appropriate symmetry operation.
One consequence of this feature is that many properties of the symmetric
Jack and Macdonald polynomials can be obtained from the corresponding
properties of the non-symmetric polynomials \cite{BF99,Ma99}.

There are, however, a number of properties of the symmetric Jack and
Macdonald polynomials which have no known relation to properties of the
non-symmetric polynomials. One example is the so-called Pieri formula
\cite{St89,Ma95,Ka97g}. To present this formula requires some notation. Let
$\kappa$ and $\lambda$ be partitions described by their diagrams and
suppose $\kappa \subset \lambda$. A skew diagram $\lambda/\kappa$ is said
to be a vertical $m$-strip if it consists of $m$ boxes, all of which
are in distinct rows. For $\lambda/\kappa$ a vertical $m$-strip define 
$\chi_m$ by $\lambda = \kappa + \chi_m$, and put
$$
U^{(\alpha)}(\lambda/\kappa) := \frac{f_N^1(\alpha \kappa + \chi_m)
f_N^{1/\alpha}(\kappa)}{f_N^1(\alpha \kappa) f_N^{1/\alpha}
(\kappa + \chi_m)}
$$
where
$$
f_n^r(\kappa) = \prod_{1 \le i < j \le n}
{((j-i)r + \kappa_i - \kappa_j)_r \over ((j-i)r)_r}, \qquad
(u)_r := {\Gamma (u+r) \over \Gamma (u)}.
$$
With this notation the Pieri formula reads
\begin{equation}\label{pf1}
e_p(z) P_\kappa(z) = \sum_{\lambda \atop \lambda/\kappa \: {\rm
a \: vertical} \: m{\rm -strip}} U^{(\alpha)}(\lambda/\kappa)
P_\lambda(z)
\end{equation}
where
$$
e_p(z) := \sum_{1 \le i_1 < \cdots < i_p \le N}
z_{i_1} \cdots z_{i_p}
$$
denotes the $p$th elementary symmetric function, and $P_\kappa(x) :=
P_\kappa(x;\alpha)$ denotes the symmetric Jack polynomial indexed by
the partition $\kappa$ and normalized so that 
when expanded in terms of monomial symmetric functions
the coefficient of the
monomial symmetric function $m_\kappa$ is unity.

It is the objective of this paper to investigate non-symmetric
analogues of the Pieri formula (\ref{pf1}). Our original idea was to
adapt the method used by Knop and Sahi \cite{KS96} to derive (\ref{pf1}),
which involves the theory of the so-called shifted Jack polynomials. 
This was passed on to D.~Marshall, who subsequently \cite{Ma00} obtained
the explicit form of the coefficients in the expansions
\begin{eqnarray}
z_i E_\eta(z) & = & \sum_{\nu: |\nu| = |\eta| + 1}
c_{\eta \nu}^{(i)} E_\nu(z) \label{m1}\\
\Big ( \sum_{i=1}^N z_i \Big ) E_\eta(z) & = & \sum_{\nu: |\nu| = |\eta| + 1}
C_{\eta \nu} E_\nu(z) . \label{m2}
\end{eqnarray}
In this work we will give an inductive proof of the evaluation of the
$c_{\eta \nu}^{(i)}$ 
which avoids all reference to the theory of the shifted Jack polynomials
(the evaluation of the $C_{\eta \nu}$ follows as a simple corollary
from knowledge of the $c_{\eta \nu}^{(i)}$).

In Section 3 of the paper we present necessary conditions on $\nu$ for the
coefficients in the expansion
\begin{equation}\label{m3}
z_{i_1} \cdots z_{i_p} E_\eta(z) =
\sum_{\nu: |\nu| = |\eta| + p} c_{\eta \nu}^{(i_1, \dots, i_p)} E_\nu(z)
\end{equation}
to be non-zero. Here use is made of the theory of shifted Jack polynomials.
In Section 4 the result of Marshall for the explicit
value of $c_{\eta \nu}^{(i)}$ is revised, and in Section 5 we present our
inductive proof of this result. The expansion (\ref{m3}) in the case
$p = N-1$, where $N$ is the number of variables $z:= (z_1,\dots,z_N)$,
is given in Section 6. In the final section, Section 7, a coefficient
in the expansion of $e_p(z) E_\eta(z)$ as a series in
$\{E_\nu\}$ is evaluated for a special value of $\nu$ and the form
of the evaluation further explored for a larger class of $\nu$.

\section{The non-symmetric Jack polynomials}
\setcounter{equation}{0}
The non-symmetric Jack polynomials $E_\eta(z)$ can be specified as 
the simultaneous polynomial eigenfunctions of the commuting operators
$$
\xi_i := \alpha z_i {\partial \over \partial z_i} +
\sum_{p<i} {z_i \over z_i - z_p} (1 - s_{ip}) +
\sum_{p > i} {z_p \over z_i - z_p} (1 - s_{ip}) + 1 - i,
$$
where $s_{ip}$ is the operator which permutes $z_i$ and $z_p$,
satisfying the eigenvalue equations
\begin{equation}\label{E1}
\xi_i E_\eta  =  \bar{\eta}_i E_\eta, \quad (i=1,\dots,N)
\end{equation}
and with coefficient of $z^\eta = z^{\eta_1} \cdots z^{\eta_N}$ unity. 
For a given composition $\eta := (\eta_1,\dots, \eta_N)$, the eigenvalue
$\bar{\eta}_i$ in (\ref{E1}) is given by
\begin{equation}\label{ac1}
\bar{\eta}_i  :=  \alpha \eta_i - \# \{ k < i | \eta_k \ge \eta_i \}
- \# \{ k > i| \eta_k > \eta_i \}.
\end{equation}

An alternative characterization of the non-symmetric Jack polynomials is
as multivariable orthogonal polynomials. 
With $z_j := e^{2 \pi i x_j}$, introduce the inner product
\begin{equation}\label{fg}
\ml f | g \mg := \int_0^1 dx_1 \cdots \int_0^1 dx_N \,
\prod_{1 \le j < k \le N} |z_k - z_j |^{2/\alpha}
f^*(z_1,\dots,z_N) g(z_1,\dots,z_N),
\end{equation}
where the $*$ denotes complex conjugation.
Suppose $|\eta| = |\nu|$ for compositions $\eta \ne \nu$.
Introduce the dominance partial ordering $<$ on
compositions by the statement that $\nu < \eta$ if
$\sum_{j=1}^p \nu_j < \sum_{j=1}^p \eta_j $
for each $p=1,\dots,N$.
Let $\eta^+$ denote the partition corresponding to the composition
$\eta$. Introduce a 
further partial ordering $\triangleleft$ by the statement that
$\nu \triangleleft \eta$ if $\nu^+ < \eta^+$, or in the case $\nu^+ = \eta^+$,
if $\nu < \eta$. Then for a given value of $|\eta|$, the $E_\eta$ can
be constructed via a Gram-Schmidt procedure from the requirements that
\begin{equation}\label{fg1}
\ml E_\eta | E_\nu \mg = 0,
\end{equation}
for $\eta \ne \nu$, and that
\begin{equation}\label{2.8'}
E_\eta(z) = z^\eta + \sum_{\nu \triangleleft \eta} c_{\eta \nu } z^\nu.
\end{equation}

We will have future use for the explicit value of
$$
{\cal N}_\eta := \ml E_\eta | E_\eta \mg.
$$
This requires the introduction of further quantities for its presentation.
Following \cite{Sa96}, define the arm and leg lengths at the node
$(i,j)$ of the diagram of a composition $\eta$ by
\begin{equation}\label{leg}
a(i,j) = \eta_i - j, \qquad
l(i,j) = \# \{ k < i | j \le \eta_k + 1 \le \eta_i \} +
\# \{ k > i | j \le \eta_k \le \eta_i \}
\end{equation}
and put
\begin{equation}\label{d1}
d_\eta' := \prod_{(i,j) \in \eta} \Big ( \alpha( a(i,j) + 1 ) +
l (i,j) \Big ), \qquad
d_\eta := \prod_{(i,j) \in \eta} \Big ( \alpha( a(i,j) + 1 ) +
l (i,j) +1 \Big ).
\end{equation}
Also, define the generalized factorial by
$$
[u]_{\eta^+}^{(\alpha)} =
\prod_{j=1}^N {\Gamma(u - (j-1)/\alpha + \eta_j^+) \over
\Gamma(u - (j-1)/\alpha)}
$$ and put
\begin{equation}\label{e1}
e_\eta = \alpha^{|\eta|} [1 + N/\alpha]_{\eta^+}^{(\alpha)}, \qquad
e_\eta' = \alpha^{|\eta|} [1 + (N-1)/\alpha]_{\eta^+}^{(\alpha)}.
\end{equation}
In terms of the quantities (\ref{d1}) and (\ref{e1}) we have
\cite{Op95,BF99}
\begin{equation}\label{en}
{{\cal N}_\eta \over {\cal N}_{(0^N)}} =
{d_\eta' e_\eta \over d_\eta e_\eta'}.
\end{equation}

Starting with $E_{(0^N)}(z) = 1$, the non-symmetric Jack polynomials can be
recursively generated from the action of just two fundamental operators.
The first of these operators is the elementary permutation operator
$s_i := s_{i \, i+1}$, which permutes $z_i$ and $z_{i+1}$. It has the
action \cite{Op95}
\begin{equation}\label{u1}
s_i E_\eta(z) = \left\{ \begin{array}{ll}
{1 \over \bar{\delta}_{i,\eta}} E_\eta(z) + \Big ( 1 - {1 \over
\bar{\delta}_{i,\eta}^2} \Big ) E_{s_i \eta}(z), & \eta_i > \eta_{i+1} \\[.2cm]
E_\eta(z), & \eta_i = \eta_{i+1} \\[.2cm]
{1 \over \bar{\delta}_{i,\eta}} E_\eta(z) + 
E_{s_i \eta}(z), & \eta_i < \eta_{i+1}
\end{array} \right.
\end{equation}
where
\begin{equation}\label{u1a}
\bar{\delta}_{i,\eta} := \bar{\eta}_i -  \bar{\eta}_{i+1}.
\end{equation}
The second required operator is the raising type operator, defined when acting
on functions according to 
$$
\Phi f(z_1,\dots, z_N) = z_N f(z_N,z_1,\dots, z_{N-1}),
$$
which has the property \cite{KS97}
\begin{equation}\label{u2}
\Phi E_\eta (z) = E_{\Phi \eta}(z), \qquad \Phi \eta := (\eta_2, \dots,
\eta_N, \eta_1 + 1).
\end{equation}
Starting from $\eta = (0^N)$, all compositions can be generated by the
action of $\Phi \eta$ and $s_i \eta$, so (\ref{u1}) and (\ref{u2})
provide the recursive generation of all the $E_\eta$.

Future use will be made of the quantities (\ref{d1}) and (\ref{e1})
with $\eta$ replaced by $s_i \eta$ and $\Phi \eta$. In particular, we require
the formulas \cite{Sa96}
$$
e_{s_i \eta} = e_\eta, \quad e_{s_i \eta}' = e_\eta', \quad
{d_{s_i \eta} \over d_\eta} = \left \{
\begin{array}{ll}
{\bar{\delta}_{i,\eta} + 1 \over \bar{\delta}_{i,\eta}}, & \eta_i > \eta_{i+1}
\\[.2cm]
{ \bar{\delta}_{i,\eta} \over
\bar{\delta}_{i,\eta} - 1 }, &
 \eta_i < \eta_{i+1} \end{array} \quad
{d_{s_i \eta}' \over d_\eta'}\right. = \left \{
\begin{array}{ll}
{\bar{\delta}_{i,\eta} \over \bar{\delta}_{i,\eta}-1}, & \eta_i > \eta_{i+1}
\\[.2cm] { \bar{\delta}_{i,\eta}+1 \over
\bar{\delta}_{i,\eta}}, &
 \eta_i < \eta_{i+1} \end{array} \right.
$$
\begin{equation}\label{15.red}
{d_{\Phi \eta} \over d_\eta} = {e_{\Phi \eta} \over e_\eta} =
\bar{\eta}_1 + \alpha + N, \qquad
{d_{\Phi \eta}' \over d_\eta'} = {e_{\Phi \eta}' \over e_\eta'}
= \bar{\eta}_1 + \alpha + N - 1.
\end{equation}

Let us now revise some aspects of the theory of non-symmetric shifted 
Jack polynomials $E_\eta^*$ \cite{Kn97}. The polynomial $E_\eta^*(z)$ is
the unique polynomial of degree $\le |\eta|$ with the property
$$
E_\eta^*(\bar{\rho}/\alpha) = 0, \qquad |\rho| \le |\eta|, \: \: \rho
\ne \eta
$$
and $E_\eta^*(\bar{\eta}/\alpha) \ne 0$ with coefficient of $z^\eta$ 
in its monomial expansion unity
($\bar{\eta} := (\bar{\eta}_1,\dots,\bar{\eta}_N)$ where the $\bar{\eta}_j$
are specified by (\ref{ac1})). The non-symmetric Jack polynomial 
$E_\eta$ is the leading homogeneous term of $E_\eta^*$ so that
\begin{equation}\label{lh}
E_\eta^*(z) = E_\eta(z) + {\rm lower \: degree \: terms}.
\end{equation}

A fundamental property of the $E_\eta^*$ is the extra vanishing condition.
Introduce the partial ordering $\preceq$ on compositions by writing
$\nu \preceq \eta$ if there exisits a permutation $\pi$ such that
$\nu_i < \eta_{\pi(i)}$ for $i < \pi(i)$ and
$\nu_i \le \eta_{\pi(i)}$ for $i \ge \pi(i)$. Note that for $\nu$ and
$\eta$ partitions the statement $\nu \preceq \eta$ is equivalent to
$\nu \subseteq \eta$ (inclusion of diagrams) but for compositions,
although $\nu \subseteq \eta$ implies $\nu \preceq \eta$ (take $\pi$ to be the
identity), the converse is not true in general. The extra vanishing condition
states \cite{Kn97}
\begin{equation}\label{v0}
E_\eta^*(\bar{\nu}/\alpha) = 0 \quad {\rm for} \quad \eta \not\prec \eta.
\end{equation}

\section{Structure of the Pieri type expansions for the non-symmetric
Jack polynomials}
\setcounter{equation}{0}
Our interest is in the coefficients $c_{\eta \nu}^{(i_1, \dots, i_p)}$
in the expansion (\ref{m3}). In this section we will use the theory of the
non-symmetric shifted Jack polynomials to present necessary conditions
for the coefficients to be non-zero.

Now the extra vanishing condition (\ref{v0}) implies that any analytic
function vanishing on 
$\{\bar{\rho}/\alpha: \eta \not\prec \rho \}$ can be written
in the form
\begin{equation}\label{v1}
f(z) = \sum_{\nu: \eta \preceq \nu} c_{\eta \nu} E_\nu^*(z).
\end{equation}
It follows from this that
\begin{equation}\label{v4}
z_{i_1} \cdots z_{i_p} E_\eta^*(z) = \sum_{\nu: \eta \preceq \nu
\atop |\nu| \le |\eta| + p} c_{\eta \nu}^{(i_1, \dots, i_p)}
E_\nu^*(z)
\end{equation}
for some coefficients $c_{\eta \nu}^{(i_1, \dots, i_p)}$. Taking the leading
homogeneous term on both sides using (\ref{lh}) gives
\begin{equation}\label{v4'}
z_{i_1} \cdots z_{i_p} E_\eta(z) = \sum_{\nu: \eta \preceq \nu
\atop |\nu| = |\eta| + p} c_{\eta \nu}^{(i_1, \dots, i_p)}
E_\nu(z)
\end{equation}
which is a refinement of (\ref{m3}).

The statement (\ref{v4'}) can be further refined by making use of 
the orthogonality (\ref{fg1}). Applying this orthogonality in (\ref{v4'})
shows that
\begin{equation}\label{f1}
c_{\eta \nu}^{(i_1, \dots, i_p)} =
{\ml E_\nu | z_{i_1} \cdots z_{i_p} E_\eta \mg \over \ml E_\nu |
E_\nu \mg}.
\end{equation}
Using the facts that 
with $z^1 := z_1 \cdots z_N$ we have
$E_{\eta + (1^N)}(z) = z^1 E_\eta(z)$,
$\ml z^1 f| z^1 g \mg = \ml f|g \mg$ and $\ml f | g \mg = \ml g | f \mg$
(the latter provided $f$ and $g$ have real coefficients) it follows from
(\ref{f1}) that
\begin{equation}\label{f2}
c_{\eta \nu}^{(i_1, \dots, i_p)} =
{\ml z_{j_1} \cdots z_{j_{N-p}} E_\nu | E_{\eta+(1^N)} \mg \over
\ml E_\nu | E_\nu \mg} =
{\ml E_{\eta+(1^N)} | z_{j_1} \cdots z_{j_{N-p}} E_\nu \mg \over
\ml E_\nu | E_\nu \mg} =
c_{\nu \eta+(1^N)}^{(j_1,\dots, j_{N-p})} {\ml E_\eta | E_\eta \mg \over
\ml E_\nu | E_\nu \mg}
\end{equation}
where $j_1,\dots,j_{N-p}$ are such that $\{1,\dots,N\} =
\{i_1,\dots,i_p\} \cup \{j_1,\dots, j_{N-p} \}$. But 
according to (\ref{v4'})
$c_{\nu \, \eta+(1^N)}^{(j_1,\dots, j_{N-p})} = 0$ for 
$\nu \not\preceq \eta +(1^N)$
and thus (\ref{f2}) implies
\begin{equation}\label{f2'}
c_{\eta \nu}^{(i_1,\dots,i_p)} = 0 \quad {\rm for} \quad
\nu \not\preceq \eta + (1^N).
\end{equation}
Hence in (\ref{v4'}) we can make the additional restriction
$\nu \preceq \eta + (1^N)$, and so obtain
\begin{equation}\label{f4}
z_{i_1} \cdots z_{i_p} E_\eta(z) =
\sum_{\nu \in {\mathbb J}_{N,p}} c_{\eta \nu}^{(i_1,\dots,i_p)}
E_\nu(z)
\end{equation}
where
\begin{equation}\label{f3}
{\mathbb J}_{N,p} := \{ \nu: \eta \preceq \nu \preceq \eta + (1^N),
\: |\nu| = |\eta| + p \}.
\end{equation}
Note that by performing the sum $1 \le i_1 < \cdots < i_p \le N$ in
(\ref{f4}) we obtain
\begin{equation}\label{f4'}
e_p(z) E_\eta(z) = \sum_{\nu \in {\mathbb J}_{N,p}}
A_{\eta \nu}^{(p)} E_\nu(z)
\end{equation}
for some constants $A_{\eta \nu}^{(p)}$.

Next we seek a more explicit description of the set ${\mathbb J}_{N,p}$.
Let $w_\eta$ be the shortest element of $S_N$ (the permutations of
$\{1,\dots,N\}$) such that $w_\eta^{-1}(\eta)$ is a partition and 
similarly define $w_\nu$. It is straightforward to show \cite{Kn97} that
if $\nu \preceq \eta$ then the permutation $\pi$ in the definition
of the partial order can be represented $\pi = w_\nu \circ
w_\eta^{-1} =: \pi_{\nu,\eta}$. Now, members $\nu$ of the set
$\mathbb J_{N,p}$ require both $\eta \preceq
\nu$ and $\nu \preceq \eta + (1^N)$ with $|\nu| = |\eta| + p$.
For the former ordering constraint the relevant permutation is
$\pi_{\eta,\nu} = \pi_{\nu,\eta}^{-1}$. Replacing $\pi$ by
$\pi^{-1}$ in the definition of $\preceq$ shows we require
\begin{equation}\label{c1}
\eta_i < \nu_{\pi_{\nu,\eta}(i)} \quad {\rm for} \quad
i < \pi_{\nu,\eta}(i), \qquad
\eta_i \le \nu_{\pi_{\nu,\eta}(i)} \quad {\rm for} \quad
i \ge \pi_{\nu,\eta}(i).
\end{equation}
For the latter ordering constraint the relevant permutation is
$\pi_{\nu,\eta}$. Replacing $\pi$ by $\pi^{-1}$ and $i$ by $\pi(i)$ in the
definition of $\preceq$ shows we require
\begin{equation}\label{c2}
\nu_{\pi_{\nu,\eta}(i)}
< \eta_i  + 1  \quad {\rm for} \quad
\pi_{\nu,\eta}(i) < i, \qquad
\nu_{\pi_{\nu,\eta}(i)} \le \eta_i + 1 \quad {\rm for} \quad
\pi_{\nu,\eta}(i) \ge i.
\end{equation}

Combining (\ref{c1}) and (\ref{c2}) gives
$$
\eta_i < \nu_{\pi_{\nu,\eta}(i)} \le \eta_i + 1
\quad {\rm for} \quad i < \pi_{\nu,\eta}(i), \qquad
\eta_i \le \nu_{\pi_{\nu,\eta}(i)} <  \eta_i + 1
\quad {\rm for} \quad i > \pi_{\nu,\eta}(i),
$$
and so
\begin{equation}\label{c3}
\nu_{\pi_{\nu,\eta}(i)} = \eta_i + 1 \quad {\rm for} \quad i < \pi_{\nu,\eta}(i), \qquad
\nu_{\pi_{\nu,\eta}(i)} = \eta_i \quad {\rm for} \quad i >
\pi_{\nu,\eta}(i).
\end{equation}
In the case $i = \pi_{\nu,\eta}(i)$ (\ref{c1}) and (\ref{c2}) give
$\eta_i \le \nu_{\pi_{\nu,\eta}(i)} \le \eta_i + 1$ and so
\begin{equation}\label{c4}
\nu_{\pi_{\nu,\eta}(i)} = \eta_i \quad {\rm or} \quad
\nu_{\pi_{\nu,\eta}(i)} = \eta_i + 1.
\end{equation}

\begin{figure}
\epsfxsize=10cm
\centerline{\epsfbox{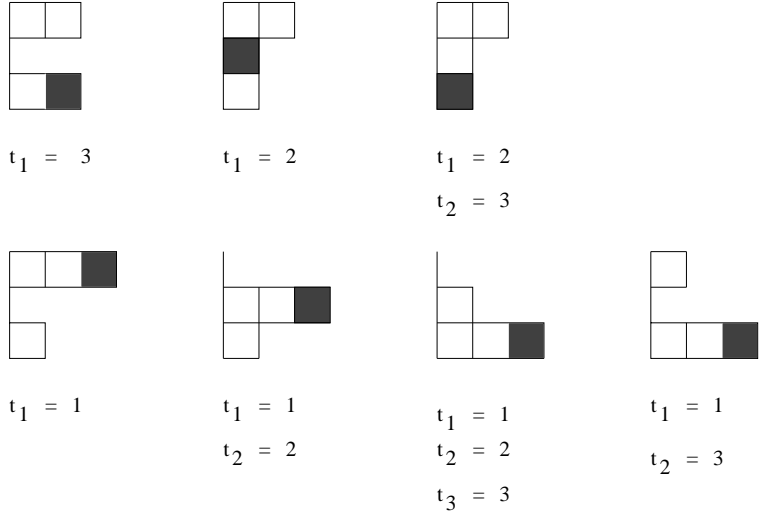}}
\caption{\label{fig1.f1} Construction of the composition
$\nu: \, \eta \preceq \nu \preceq \eta + (1)^3$ with $\eta
= (2,0,1)$ and $|\nu| = |\eta| + 1$. The unshaded boxes originate from
the diagram of $\eta$. With reference to the original diagram of
$\eta$ the row with the additional box (shaded) must
move downwards or stay stationary, while the rows with no box added
move upwards or stay stationary. The labels in the description
(\ref{c8}) are also noted.} 
\end{figure}

It remains to implement the requirement $|\nu| = |\eta| + p$. We see from
(\ref{c3}) and (\ref{c4}) that we must have
\begin{equation}\label{c5}
\nu_{\pi_{\nu,\eta}(i_r)} = \eta_{i_r} + 1 \qquad (r=1,\dots,p)
\end{equation}
for some $1 \le i_1 < \cdots < i_p \le N$ and
\begin{equation}\label{c6}
\nu_{\pi_{\nu,\eta}(j_r)} = \eta_{j_r} \qquad (r=1,\dots,N-p)
\end{equation}
where $\{i_1,\dots,i_p\} \cup \{j_1,\dots,j_{N-p} \} =
\{1,2,\dots, N\}$. Combining (\ref{c5}) and (\ref{c6}) with (\ref{c3})
and (\ref{c4}) shows compositions $\nu \in \mathbb J_{n,p}$ are
characterized by the properties
\begin{eqnarray}\label{c7}
\nu_{\pi(i_r)} = \eta_{i_r} + 1 \quad {\rm for} \quad
i_r \le \pi(i_r) \qquad r=1,\dots,p \nonumber \\
\nu_{\pi(j_r)} = \eta_{j_r} \quad {\rm for} \quad
j_r \ge \pi(j_r) \qquad r=1,\dots,N-p
\end{eqnarray}
for some permutation $\pi$ ($\pi = \pi_{\nu,\eta}$ suffices).
The characterization (\ref{c7}) can be interpreted in terms of the diagram
of $\eta$. We begin by adding one box to the rows $i_1,\dots,i_p$. Then
we consider all rearrangements of the rows such the rows with a box
added move downwards or stay stationary, while the rows with no box added
move upwards or stay stationary. An example is given in Figure
\ref{fig1.f1}.

\vspace{.5cm}

\begin{figure}
\epsfxsize=6cm
\centerline{\epsfbox{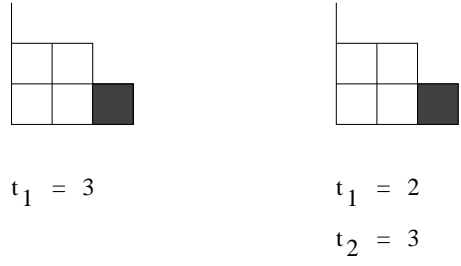}}
\caption{\label{fig1.f2} In this example, starting with
$\eta = (0,2,2)$, two different choices of subsets
$I = \{t_1,\dots, t_s\}$ give the same composition, but only the second
subset is maximal (note that in the first diagram $\eta_2 = \eta_{t_1}$).
}
\end{figure}

In the case $p=1$ the compositions $\nu$ defined by (\ref{c7}) and
thus belonging to the set $\mathbb J_{N,1}$ have the property of being
the minimal elements lying above $\eta$ \cite{Kn97}. Note that this
set can be indexed by subsets $I = \{t_1,\dots, t_s\}$ of
$\{1,\dots,N\}$ with $t_1 < \cdots < t_s$ which correspond to the element
\begin{equation}\label{c80}
\nu =: c_I(\eta) \in \mathbb J_{N,1}
\end{equation}
where
\begin{eqnarray}\label{c8}
\nu_{t_j} & = & \eta_{t_{j+1}} \qquad j=1,\dots, s-1 \nonumber \\
\nu_{t_s} & = & \eta_{t_1} + 1  \nonumber \\
\nu_i & = & \eta_i \qquad i \notin I
\end{eqnarray}
Furthermore, the subset $I$ is called maximal with respect to $\eta$ if
$I \neq \emptyset$ and
\begin{eqnarray}\label{c8a}
\eta_j & \ne & \eta_{t_u} \qquad j=t_{u-1} + 1, \dots, t_u - 1
\: \: (u=1,\dots,s; \, t_0 := 0) \nonumber \\
\eta_j & \ne & \eta_{t_1} + 1 \qquad j=t_s+1, \dots, N
\end{eqnarray}
It follows from (\ref{c8}) that an equivalent way to characterize the maximal
subsets is via the conditions
\begin{eqnarray}\label{c9}
\nu_j &\ne & \nu_{t_s} - 1, \qquad j=1,\dots, t_1 - 1 \nonumber \\
\nu_j & \ne & \nu_{t_u}, \qquad j=t_u+1, \dots, t_{u+1} - 1 \: \:
(u=1,\dots,s; \: t_{s+1} := N+1).
\end{eqnarray}
It is shown in \cite{Kn97} that it is only these maximal
subsets which give distinct compositions $\nu$ (we illustrate this point
in Figure \ref{fig1.f2}).  Thus we can write
\begin{equation}\label{c8b}
\mathbb J_{N,1} :=
\mathbb J_{N,1}[\eta] =
 \{ \nu: \, \nu = c_I(\eta), \: I \: {\rm maximal} \}.
\end{equation}
It is also convenient to introduce the set $\mathbb J_\eta$ of maximal
subsets
\begin{equation}\label{c8c}
\mathbb J_\eta = \{I: \, I \: {\rm is \: maximal \: w.r.t. \: } \eta \},
\end{equation}
so that $\mathbb J_{N,1} = \{c_I(\eta): \, I \in  \mathbb J_\eta\}$.

\vspace{.5cm}

\section{A Pieri type formula for the non-symmetric Jack polynomials
in the case $p=1$}
\setcounter{equation}{0}
So far the theory of shifted Jack polynomials \cite{Kn97} has been used
to deduce the structural formula (\ref{f4}), and also notions from
that theory are used to label the set $\mathbb J_{N,p}$ appearing
in (\ref{f4}) in terms of certain maximal subsets $I$. To now evaluate
the coefficients in (\ref{f4}), the most natural way to proceed is to
make further use of theory from \cite{Kn97}. In the case $p=1$
this part of the program
has recently been successfully undertaken by Marshall \cite{Ma00}.
The presentation of the result requires some notation.

First write
\begin{equation}
a(x,y) := {1 \over \alpha (x-y)}, \qquad
b(x,y) := {x-y-1/\alpha \over x - y}.
\end{equation}
For $I = \{t_1,\dots,t_s\} \subseteq \{1,\dots,N\}$, $I \ne \emptyset$,
$t_1 < \cdots < t_s$ put
\begin{eqnarray}
A_I(x) &:=& \Big ( \prod_{u=1}^{s-1} a(x_{t_u}, x_{t_{u+1}}) \Big )
a(x_{t_s} - 1, x_{t_1}) \label{a2} \\
B_I(x) &:=& \Big ( \prod_{u=1}^{s} \prod_{j=t_u+1}^{t_{u+1}-1}
b(x_{t_u},x_j) \Big ) 
(x_{t_s} + (N-1)/\alpha) \prod_{j=1}^{t_1 - 1}
b(x_{t_s} -1, x_j), \quad t_{s+1} := N + 1 \label{a3} \\
\tilde{B}_I(x) &:=& 
\Big ( \prod_{u=1}^{s} \prod_{j=t_{u-1}+1}^{t_{u}-1}
b(x_{t_u},x_j) \Big ) \Big ( \prod_{j=t_s +1}^N b(x_{t_1}+1,x_j) \Big )
(x_{t_1} +1 + (N-1)/\alpha) , \quad t_0 := 0
\label{a3'}
\end{eqnarray}
and for $i \in I$ write
\begin{eqnarray}
\chi_I^{(i)}(x) & = & \left \{ \begin{array}{ll}
\alpha(x_{t_{k-1}} - x_i), & i=t_k \: \: (k=2,\dots,s) \\
\alpha(x_{t_s} - x_i - 1), & i=t_1 \end{array} \right. 
\nonumber \label{a3a}\\
\tilde{\chi}_I^{(i)}(x) &=& \left \{
\begin{array}{ll} \alpha (x_i - x_{t_{k+1}}), & i=t_k \: \:
(k=1,\dots,s-1)  \\ 
\alpha (x_i - x_{t_1} - 1), & i = t_s. \end{array} \right. \label{a3b}
\end{eqnarray}

In terms of 
these quantities,
and the quantity
$d_\eta'$ of (\ref{d1}), the result of Marshall \cite{Ma00} reads
\begin{equation}\label{fu1}
z_i E_\eta(z) = \alpha d_\eta' \sum_{I \in \mathbb J_{\eta} \atop
{\rm given} \, i \in I}
{\chi_I^{(i)}(\overline{c_I(\eta)}/\alpha) 
A_I(\overline{c_I(\eta)}/\alpha)
B_I(\overline{c_I(\eta)}/\alpha)
\over d_{c_I(\eta)}' }
E_{c_I(\eta)}(z).
\end{equation}
Also, noting from (\ref{a3b}) that
\begin{equation}\label{4.47'}
\sum_{i \in I} \chi_I^{(i)}(x) = - \alpha
\end{equation}
it follows from (\ref{fu1}) that  \cite{Ma00}
\begin{equation}\label{fu2}
\Big ( \sum_{i=1}^N z_i \Big ) E_\eta (z) = - \alpha^2
d_\eta'
\sum_{I \in \mathbb J_\eta}
{A_I((\overline{c_I(\eta)}/\alpha) B_I((\overline{c_I(\eta)}/\alpha) \over
d_{c_I(\eta)}' } E_{c_I(\eta)}(z).
\end{equation}
We remark that it follows from the definition (\ref{a3}) of $B_I(x)$ that
for $I$ not maximal (i.e.~cases for which the relations (\ref{c9}) are
not obeyed), $B_I(\overline{c_I(\eta)}/\alpha) = 0$. Thus the restriction to
maximal subsets in the summation of (\ref{fu1}) and (\ref{fu2}) is in
fact a feature of the analytic form of the coefficients.

The dependence on $\overline{c_I(\eta)}$ in (\ref{fu1}) and (\ref{fu2}) can
be replaced by a dependence on $\bar{\eta}$. Thus we note from the
definitions (\ref{a2}), (\ref{a3'}) and (\ref{a3b}) that
$$
A_I(\overline{c_I(\eta)}/\alpha) = A_I(\bar{\eta}/\alpha), \quad
B_I(\overline{c_I(\eta)}/\alpha) = \tilde{B}_I(\bar{\eta}/\alpha), \quad
\chi_I^{(i)}(\overline{c_I(\eta)}/\alpha) =
\tilde{\chi}_I^{(i)}(\bar{\eta}/\alpha) \: \: (i \in I),
$$
which when substituted in (\ref{fu1})  and (\ref{fu2})
give
\begin{eqnarray}
z_i E_\eta(z) & = & \alpha d_\eta' \sum_{I \in \mathbb J_{\eta} \atop
{\rm given} \, i \in I}
{\tilde{\chi}_I^{(i)}(\bar{\eta}/\alpha)
A_I(\bar{\eta}/\alpha) \tilde{B}_I(\bar{\eta}/\alpha) \over d_{c_I(\eta)}' }
E_{c_I(\eta)}(z) \label{fz1} \\
\Big ( \sum_{i=1}^N z_i \Big ) E_\eta (z) & = & - \alpha^2
d_\eta'
\sum_{I \in \mathbb J_\eta}
{A_I(\bar{\eta}/\alpha) \tilde{B}_I(\bar{\eta}/\alpha) \over d_{c_I(\eta)}' }
E_{c_I(\eta)}(z). \label{fz2}
\end{eqnarray}
A still more useful form of (\ref{fz1}) results by introducing
\begin{equation}\label{a3'h}
\hat{B}_I(\bar{\eta}/\alpha)  :=  \alpha {e_\eta' \over
e_{c_I(\eta)}'} \tilde{B}_I(\bar{\eta}/\alpha) 
 =  \Big ( \prod_{u=1}^{s} \prod_{j=t_{u-1}+1}^{t_{u}-1}
b(x_{t_u},x_j) \Big ) \Big ( \prod_{j=t_s+1}^Nb(x_{t_s}+1,x_j) \Big ),
\quad t_0 := 0
\end{equation}
where the equality follows from (\ref{e1}) and (\ref{a3'}).
In terms of this quantity (\ref{fz1}) reads
\begin{equation}\label{obt}
z_i E_\eta(z) = {d_\eta' \over e_\eta'}
\sum_{I \in \mathbb J_{\eta} \atop
{\rm given} \, i \in I}
{e_{c_I(\eta)}'\tilde{\chi}_I^{(i)}(\bar{\eta}/\alpha)
A_I(\bar{\eta}/\alpha) \hat{B}_I(\bar{\eta}/\alpha) \over d_{c_I(\eta)}' }
E_{c_I(\eta)}(z).
\end{equation}

\section{Inductive proof}
\setcounter{equation}{0}
In this section we will provide an inductive proof of (\ref{obt}). This
has the advantage of being independent of the theory of the shifted
Jack polynomials, relying only on the recurrence properties
(\ref{u1}) and (\ref{u2}) of the non-symmetric Jack polynomials
themselves.

\subsection*{Strategy}
It has already been remarked that starting with $E_{(0^N)}(z) = 1$,
the non-symmetric Jack polynomials can be generated recursively from
the recurrence properties (\ref{u1}) and (\ref{u2}). To make use of
these properties, suppose for a given $\eta$ we know the coefficients
$c_{\eta, \nu}^{(j)}$ in the expansion
\begin{equation}\label{k1}
z_j E_\eta(z) = \sum_{\nu \in \mathbb J_{N,1}[\eta]}
c_{\eta, \nu}^{(j)} E_\nu(z)
\end{equation}
for each $j=1,\dots,N$. Then, with $z_{N+1} := z_1$ and
$c_{\eta, \nu}^{(N+1)} := c_{\eta, \nu}^{(1)}$, (\ref{u2}) gives
\begin{equation}\label{k2}
z_j E_{\Phi \eta}(z) = z_j \Phi E_\eta(z) = \Phi(z_{j+1}
E_\eta(z) ) =
\sum_{\nu \in \mathbb J_{N,1}[\eta]} c_{\eta, \nu}^{(j+1)}
E_{\Phi \nu}(z).
\end{equation}
This shows $z_j E_{\Phi \eta}(z)$ can be computed from knowledge of
the expansion (\ref{k1}) for the given $\eta$. Moreover, we can can
give an explicit relationship between coefficients.
To demonstrate this,
for $I \subseteq \{1,\dots,N\}$, $I \ne \emptyset$ put
\begin{equation}\label{thi}
\Phi (I) := \{ j-1 | j \in I \cap \{2,\dots,N\} \} \cup
\{N | 1 \in I \}.
\end{equation}
Then we can check that $\Phi(I)$ is maximal with respect to $\Phi \eta$ if
and only if $I$ is maximal with respect to $\eta$. This means that in
(\ref{k2}) we can replace the summation $\nu \in \mathbb J_{N,1}[\eta]$
by $\Phi \nu \in \mathbb J_{N,1}[\Phi \eta]$, which allows us
to change variables
$\Phi \nu \mapsto \nu$ to obtain
\begin{equation}\label{k3}
z_j E_{\Phi \eta}(z) =
\sum_{\nu \in \mathbb J_{N,1}[\Phi \eta]} c_{\eta, \Phi^{-1}\nu}^{(j+1)}
E_{\Phi \nu}(z).
\end{equation}
On the other hand (\ref{k1}) gives
\begin{equation}\label{k4}
z_j E_{\Phi \eta}(z) =
\sum_{\nu \in \mathbb J_{N,1}[\Phi \eta]} c_{\Phi \eta,\nu}^{(j)}
E_{\Phi \nu}(z).
\end{equation}
Comparing (\ref{k3}) and (\ref{k4}) shows we require
\begin{equation}\label{k5}
c_{\Phi \eta,\Phi\nu}^{(j)} = c_{\eta, \nu}^{(j+1)}, \quad
(j=1,\dots, N-1) \qquad
c_{\Phi \eta,\Phi\nu}^{(N)} = c_{\eta, \nu}^{(1)}.
\end{equation}

Let us now consider the computation of $z_j E_{s_i \eta}(z)$ for
$\eta_i < \eta_{i+1}$ from knowledge of the expansion (\ref{k1}) for
the given $\eta$. For this purpose we rewrite (\ref{k1}) as
\begin{equation}\label{k1a}
z_j E_\eta(z) = \sum_\nu \alpha_{\eta,\nu}^{(j)} E_\nu(z)
\end{equation}
where $\alpha_{\eta,\nu}^{(j)} = c_{\eta,\nu}^{(j)}$,
$\nu \in \mathbb J_{N,1}[\eta]$ and $\alpha_{\eta,\nu}^{(j)} = 0$
otherwise. By doing this the sum over $\nu$ in (\ref{k1a}) is
unrestricted. Since from (\ref{u1}), with $\eta_i < \eta_{i+1}$,
\begin{eqnarray}\label{s1}
z_j E_{s_i \eta}(z)
 & = & z_j \Big ( s_i E_\eta(z) - {1 \over \bar{\delta}_{i,
\eta}} E_\eta(z) \Big ) \nonumber \\
&=& \left \{ \begin{array}{ll}
s_i(z_j E_\eta(z)), & j \ne i,i+1 \\
s_i(z_i E_\eta(z)), & j = i+1 \\
s_i(z_{i+1} E_\eta(z)), & j = i \end{array} \right. -
{1 \over \bar{\delta}_{i,\eta}} z_j E_\eta(z)
\end{eqnarray}
we see that knowledge of $z_j E_\eta(z)$ for each $j=1,\dots,N$ implies
the value of $z_j E_{s_i \eta}(z)$. We want to exhibit this feature as
a recurrence for the coefficients $\alpha_{\eta,\nu}^{(j)}$. Now, from
(\ref{k1a}) and (\ref{u1})
\begin{eqnarray}
s_i(z_j E_\eta(z)) & = & \sum_{\nu_i < \nu_{i+1}}
\{ \alpha_{\eta,\nu}^{(j)} \bar{\delta}_{i,\nu}^{-1} +
\alpha_{\eta,s_i\nu}^{(j)}(1 - \bar{\delta}_{i,\nu}^{-2}) \} E_\nu(z) +
\sum_{\nu_i = \nu_{i+1}} \alpha_{\eta,\nu}^{(j)} E_\nu(z) \nonumber \\
&& + \sum_{\nu_i > \nu_{i+1}} 
\{ \alpha_{\eta,\nu}^{(j)} \bar{\delta}_{i,\nu}^{-1} +
\alpha_{\eta,s_i\nu}^{(j)} \} E_\nu(z) \label{s2}
\end{eqnarray}
while (\ref{k1a}) itself gives
\begin{equation}\label{s3}
z_j E_{s_i \eta}(z) = \sum_{\nu} \alpha_{s_i \eta, \nu}^{(j)}
E_\nu(z).
\end{equation}
Substituting (\ref{s3}), (\ref{s2}) and (\ref{k1a}) in (\ref{s1}) and
equating coefficients of $E_\nu(z)$ gives a recurrence allowing 
$\alpha_{\eta,s_i\nu}^{(j)}$ to be computed. In the recurrence it is
necessary to distinguish the cases $\nu_i < \nu_{i+1}$ from
$\nu_i > \nu_{i+1}$. However this can be avoided if we write the
recurrence in terms of the
quantity
$$
\tilde{\alpha}_{\eta,\nu}^{(j)} := {d_\nu' e_\eta' \over d_\eta' e_\nu'}
\alpha_{\eta,\nu}^{(j)}
$$
and make use of (\ref{15.red}). We then find
for $\nu_i \ne \nu_{i+1}$
\begin{eqnarray}
(1 + \bar{\delta}_{i,\eta}^{-1}) \tilde{\alpha}_{s_i \eta, \nu}^{(j)} & = &
(1 - \bar{\delta}_{i,\nu}^{-1} ) \tilde{\alpha}_{\eta,s_i\nu}^{(j)}
+ (\bar{\delta}_{i,\nu}^{-1} - \bar{\delta}_{i,\eta}^{-1})
\tilde{\alpha}_{\eta,\nu}^{(j)} \quad (j \ne i,i+1 )
 \nonumber \\
(1 + \bar{\delta}_{i,\eta}^{-1}) \tilde{\alpha}_{s_i \eta, \nu}^{(i+1)} & = &
(1 - \bar{\delta}_{i,\nu}^{-1} ) \tilde{\alpha}_{\eta,s_i\nu}^{(i)} +
\bar{\delta}_{i,\nu}^{-1}\tilde{\alpha}_{\eta,\nu}^{(i)} -
\bar{\delta}_{i,\eta}^{-1} \tilde{\alpha}_{\eta,\nu}^{(i+1)} \nonumber \\
(1 + \bar{\delta}_{i,\eta}^{-1}) \tilde{\alpha}_{s_i \eta, \nu}^{(i)} & = &
(1 - \bar{\delta}_{i,\nu}^{-1} ) \tilde{\alpha}_{\eta,s_i\nu}^{(i+1)} +
\bar{\delta}_{i,\nu}^{-1}\tilde{\alpha}_{\eta,\nu}^{(i+1)} -
\bar{\delta}_{i,\eta}^{-1} \tilde{\alpha}_{\eta,\nu}^{(i)} \label{s4}
\end{eqnarray}
while for $\nu_i = \nu_{i+1}$
\begin{eqnarray}
(1 + \bar{\delta}_{i,\eta}^{-1}) \tilde{\alpha}_{s_i \eta, \nu}^{(j)} & = &
(1 - \bar{\delta}_{i,\eta}^{-1} ) \tilde{\alpha}_{\eta,\nu}^{(j)} 
\quad (j \ne i,i+1 )
 \nonumber \\
(1 + \bar{\delta}_{i,\eta}^{-1}) \tilde{\alpha}_{s_i \eta, \nu}^{(i+1)} & = &
\tilde{\alpha}_{\eta,\nu}^{(i)} - \bar{\delta}_{i,\eta}^{-1}
\tilde{\alpha}_{\eta, \nu}^{(i+1)} \nonumber \\
(1 + \bar{\delta}_{i,\eta}^{-1}) \tilde{\alpha}_{s_i \eta, \nu}^{(i)} & = &
\tilde{\alpha}_{\eta,\nu}^{(i+1)} - \bar{\delta}_{i,\eta}^{-1}
\tilde{\alpha}_{\eta, \nu}^{(i)}. \label{s5}
\end{eqnarray}
Noting that for $\nu_i = \nu_{i+1}$ we have $s_i \nu = \nu$, we see that
the equations (\ref{s4}) remain valid in that they reduce to
 the equations (\ref{s5}), 
so it suffices to consider (\ref{s4}) for all $\nu$.

Starting from knowledge of $c_{(0^N),\nu}^{(j)}$ and $\tilde{\alpha}_{(0)^N,
\nu}^{(j)}$ the recurrences (\ref{k5}) and (\ref{s4}) can be
used to compute all the $c_{\eta,\nu}^{(j)}$ and $\tilde{\alpha}_{\eta,
\nu}^{(j)}$. Thus, after independently establishing their validity in
the case $\eta = (0^N)$, we want to show the functional forms
\begin{eqnarray}
\tilde{\alpha}_{\eta,\nu}^{(j)} & = & \left \{\begin{array}{ll}
 \tilde{\chi}_I^{(j)}(\bar{\eta}/\alpha)
A_I(\bar{\eta}/\alpha) \hat{B}_I(\bar{\eta}/\alpha), &
I \in \mathbb J_\eta \: \: {\rm and} \: \: j \in I \: \: {\rm with}
\: \: \nu = c_I(\eta) \\
0, & {\rm otherwise}
\end{array} \right. \label{s5a} \\
c_{\eta, \nu}^{(j)} & = & {d_\eta' e_\nu' \over d_\nu' e_\eta'}
\tilde{\chi}_I^{(j)}(\bar{\eta}/\alpha)
A_I(\bar{\eta}/\alpha) \hat{B}_I(\bar{\eta}/\alpha), \qquad
\nu = c_I(\eta), \label{s5b'}
\end{eqnarray}
with $\tilde{\chi}^{(j)}$, $A_I$ and $\hat{B}_I$ as specified by
(\ref{a3b}), (\ref{a2}) and (\ref{a3'h}) respectively, satisfy the
recurrences (\ref{k5}) and (\ref{s4}) as appropriate.

\subsection*{Verification of the initial conditions}
We require the expansion of $z_i$ in terms of $\{E_\nu\}$. From the
recurrences (\ref{k5}) and (\ref{s4}) we can readily show
$$
E_{(0^k 1 0^{N-k-1})}(z) = z_{k+1} + {1 \over \alpha + k + 1}
(z_{k+2} + \cdots + z_N).
$$
Thus the expansion of $\{ E_{(0^k 1 0^{N-k-1})} \}$, $k=0,\dots, N-1$,
in terms of $\{z_i\}$ has a triangular structure. This makes the task of
inverting the formulas straightforward, provided we start with $z_N$
and then compute the expansion of $z_{N-1}$ etc.. We find
\begin{equation}\label{ag}
z_i = E_{(0^{i-1} 1 0^{N-i})}(z) - {1 \over \alpha + i} 
E_{(0^i 1 0^{N-i-1})}
- {1 \over \alpha + i + 1} E_{(0^{i+1} 1 0^{N-i-2})}(z) - \cdots -
{1 \over \alpha + N - 1} E_{(0^{N-1} 1)}(z).
\end{equation}

On the other hand for $\nu = c_I(0^N)$, $I \in \mathbb J_{(0^N)}$ and
$i \in I$ we see from (\ref{c8}) that the only possibilities are
$\nu = (0^{j-1} 1 0^{N-j})$ with $j \ge i$ so that (\ref{s5a}) gives
\begin{equation}\label{ag1}
z_i = \sum_{j=i}^N c_j^{(i)} E_{(0^{j-1} 1 0^{N-j})}(z)
\end{equation}
for some constants $c_j^{(i)}$. The general structure of (\ref{ag1})
is in agreement with (\ref{ag}). To check that the coefficients
agree we note from (\ref{c8}) that for $I \in
\mathbb J_{(0^N)}$, we must have $I = \{1,2,\dots,j\}$ and $t_l = l$
$(l=1,\dots,j)$. Noting also that with $\eta = (0^N)$,
$\bar{\eta}_i = - (i-1)$, we see that (\ref{a3b}), (\ref{a2}) and (\ref{obt})
give
\begin{eqnarray}\label{subs}
\tilde{\chi}_I^{(i)}(\overline{(0^N)}/\alpha) & = &
\left \{ \begin{array}{ll} 1, & i < j \\
- (N - 1 + \alpha), & i = j \end{array} \right.
\nonumber \\
A_I(\overline{(0^N)}/\alpha) & = & - {1 \over j - 1 + \alpha}, \nonumber \\ 
\hat{B}_I(\overline{(0^N)}/\alpha) & = & {j - 1 + \alpha \over N - 1 + \alpha}.
\end{eqnarray}
Substituting (\ref{subs}), together with the evaluations
$$
d_{(0^N)}' = e_{(0^N)}'=1, \quad
e_{(0^{j-1} 1 0^{N-j})}' = \alpha + N - 1, \quad
d_{(0^{j-1} 1 0^{N-j})}' = \alpha + j - 1
$$
from (\ref{e1}) and (\ref{d1}), in (\ref{s5a}) we see that the
coefficients in (\ref{ag1}) are as required by (\ref{ag}).

\subsection*{Verification of the recurrences}
Consider (\ref{k5}). From the definitions (\ref{ac1}) and (\ref{u2}) we
can check
\begin{equation}\label{ot}
(\overline{\Phi \eta})_i = (\bar{\eta})_{i+1}, \quad i \ne N \qquad
(\overline{\Phi \eta})_N = (\bar{\eta})_1 + \alpha.
\end{equation}
With $\Phi(I)$ defined by (\ref{thi}) and
$I = \{t_1,\dots,t_s\}$, making use of (\ref{ot}) it follows from the
definition (\ref{a2}) that
\begin{eqnarray}
\lefteqn{
A_{\Phi(I)}(\overline{\Phi \eta}/ \alpha)} \nonumber
\\&&  =  \left \{ \begin{array}{ll}
\Big ( \prod_{u=1}^{s-1}a((\overline{\Phi \eta})_{t_u - 1}/\alpha,
(\overline{\Phi \eta})_{t_{u+1} - 1}/\alpha) \Big ) 
a((\overline{\Phi \eta})_{t_s - 1}/\alpha -1,
(\overline{\Phi \eta})_{t_1 - 1}/\alpha), & t_1 \ne 1 \\[.2cm]
(a((\overline{\Phi \eta})_{t_{s-1}-1})/\alpha,
(\overline{\Phi \eta})_N/\alpha)
a((\overline{\Phi \eta})_N/\alpha - 1, (\overline{\Phi \eta})_{t_2-1}/
\alpha)
\prod_{u=2}^{s-1}a((\overline{\Phi \eta})_{t_{u} - 1}/\alpha,
(\overline{\Phi \eta})_{t_{u+1} - 1}/\alpha) , & t_1 =1 \end{array} \right.
\nonumber \\
&& =  \prod_{u=1}^{s-1}a(\bar{\eta}_{t_u}/\alpha,
\bar{\eta}_{t_{u+1}}/\alpha) \, a(\bar{\eta}_{t_s}/\alpha - 1,
\bar{\eta}_{t_1}/\alpha) \: = \:
A_I(\bar{\eta}/\alpha).
\end{eqnarray}
Similar calculations show
$$
\hat{B}_{\Phi(I)}(\overline{\Phi \eta}/ \alpha) =
\hat{B}_{I}(\bar{\eta}/\alpha), \quad
\tilde{\chi}_{\Phi({I})}^{(i)}(\overline{\Phi \eta}/ \alpha) =
\tilde{\chi}_{I}^{(i+1)}(\bar{\eta}/\alpha), \: (i=1,\dots,
N-1),
$$
$$
 \tilde{\chi}_{\Phi({I})}^{(N)}(\overline{\Phi \eta}/ \alpha) =
\tilde{\chi}_{I}^{(1)}(\bar{\eta}/\alpha), \quad
c_{\Phi(I)}(\Phi \eta) =
\Phi (c_I(\eta)) .
$$ 
These formulas together with the appropriate
formula from (\ref{15.red})
immediately imply (\ref{k5}) is satisfied by (\ref{s5a}).

The recurrences (\ref{s4}) are not so straightforward.
One complication is that the cases
\begin{equation}
{\rm (i) } \: i,i+1 \notin I \quad
{\rm (ii) } \: i \in I, \: i+1 \notin I \quad
{\rm (iii) } \: i \notin I, \: i+1 \in I \quad
{\rm (iv) } \: i, i+1 \in I
\label{Aiv}
\end{equation}
must be treated separately, in addition to the division of cases depending
on the value of $(j)$.
Independent of the division of cases (\ref{Aiv}), the fact that 
$\tilde{\alpha}_{\eta,\nu}^{(j)} = 0$ for $\nu \notin
\mathbb J_{N,1}[\eta]$ used in (\ref{s4}) gives
\begin{equation}\label{s4m}
\tilde{\alpha}_{s_i \eta, \nu}^{(j)} = 0, \qquad
\nu \ne c_I(\eta), \: I \in \mathbb J_\eta \: \: {\rm and} \: \:
\nu \ne s_i c_I(\eta), \: I \in \mathbb J_\eta.
\end{equation}
For  (\ref{s4m}) 
to be consistent with (\ref{s5a}) we must show
\begin{equation}\label{state} 
c_{I'}(s_i \eta) = c_I(\eta) \quad
{\rm or} \quad c_{I'}(s_i \eta) = s_i c_I(\eta) \: {\rm for \: some 
} \:
I \in  \mathbb J_{\eta}. 
\end{equation}
The validity of this statement will be verified for each of the
cases separately.

First suppose $I' \in \mathbb J_{s_i \eta} \Big |_{i,i+1 \notin I'}$.
The definitions (\ref{c8}) and (\ref{c8a}) give that this is equivalent
to the statement that $I'=I$, $I \in \mathbb J_\eta \Big |_{i,i+1 \in I}$,
and 
\begin{equation}\label{r1}
c_I(s_i \eta) = s_i c_I(\eta).
\end{equation}
Suppose next $I' \in \mathbb J_{s_i \eta} \Big |_{i+1 \in I', i \notin I'}$.
Then there are two possibilities. The first is $I' = (I \cup \{i+1\})
\backslash
\{i\}$, $I \in \mathbb J_\eta \Big |_{i \in I, i+1 \notin I}$, with
\begin{equation}\label{r2}
c_{(I \cup \{i+1\})\backslash\{i\}}(s_i \eta) = s_i c_I(\eta).
\end{equation}
The second is $I' = I \backslash\{i\}$, $I \in \mathbb J_\eta
\Big |_{i,i+1 \in I}$ with 
\begin{equation}\label{r3}
c_{I\backslash\{i\}}(s_i \eta) = c_I(\eta).
\end{equation}
In the case $I' \in \mathbb J_{s_i \eta} \Big |_{i \in I, i+1 \notin I'}$
the only possibility is $I'=(I \cup \{i\}) \backslash\{i+1\}$,
$I \in  \mathbb J_\eta \Big |_{i \notin I, i+1 \in I}$, with
\begin{equation}\label{r4}
c_{(I \cup \{i\}) \backslash\{i+1\}}(s_i \eta) = s_i c_I(\eta).
\end{equation}
The remaining case is $I' \in \mathbb J_{s_i \eta} \Big |_{i,i+1 \in I'}$.
Then we can have $I' = I \cup \{i\}$, $I \in \mathbb J_\eta
\Big |_{i+1 \in I, i \notin I}$ with
\begin{equation}\label{r5}
c_{I \cup \{i\}}(s_i \eta) = c_I(\eta).
\end{equation}
These results together verify (\ref{state}).  Thus we can restrict attention
to the cases 
\begin{equation}\label{cases}
\nu = c_I(\eta), \qquad \nu = s_i c_I(\eta), \quad (I \in \mathbb J_\eta).
\end{equation}

\noindent
{\it The case $i, i+1 \notin I$}

Because $\tilde{\alpha}_{\eta, \nu}^{(j)}$ requires $j \in I$ to be
non-zero, while we are considering the case $i, i+1 \notin I$, the
second and third equations in (\ref{s4}) give
$\tilde{\alpha}_{s_i \eta, \nu}^{(i)} = 0$ and
$\tilde{\alpha}_{s_i \eta, \nu}^{(i+1)} = 0$, which is consistent with
(\ref{s5a}). Thus we can restrict attention to the first equation in
(\ref{s4}). Also, if $\nu_i = \nu_{i+1}$ with $i, i+1 \notin I$ then
we must have $\eta_i = \eta_{i+1}$. Since in the induction procedure
it suffices to consider only the cases $\eta_{i+1} > \eta_i$ we can
suppose $\nu_i \ne \nu_{i+1}$.

Suppose $\nu = c_I(\eta)$. The assumptions that $i, i+1 \notin I$
and $\eta_i \ne \eta_{i+1}$
together with (\ref{r1}) imply there is no $I' \in \mathbb J_\eta$ such
that $c_{I'}(\eta) = s_i c_I(\eta)$ and thus $\tilde{\alpha}_{\eta, 
s_i \nu}^{(j)} = 0$. Also in this case it follows from the definition
(\ref{u1a}) that $\bar{\delta}_{i, c_I(\eta)} = \bar{\delta}_{i,\eta}$.
Substituting these formulas in the first equation of (\ref{s4}) gives
\begin{equation}\label{s4x}
\tilde{\alpha}_{s_i \eta, c_I(\eta)}^{(j)} = 0, \qquad
j \ne i, i+1, \quad I \in \mathbb J_\eta \Big |_{i,i+1 \notin I}.
\end{equation}
This is consistent with (\ref{s5a}) because (\ref{r1}) and the
surrounding sentence implies there is no $I' \in \mathbb J_{s_i \eta}$
such that $c_{I'}(s_i \eta) = c_I(\eta)$.

According to (\ref{cases}) the remaining possibility for a non-zero value is
$\nu = s_i c_I(\eta)$. From the above reasoning we know from this
choice of $\nu$, $\tilde{\alpha}_{\eta,\nu}^{(j)} = 0$, while (\ref{u1a})
gives $\bar{\delta}_{i,s_i c_I(\eta)} = - \bar{\delta}_{i,\eta}$.
Thus the first equation in (\ref{s4}) reduces to
\begin{equation}\label{s4y}
\tilde{\alpha}_{s_i \eta, s_i c_I(\eta)}^{(j)} =
\tilde{\alpha}_{\eta, c_I(\eta)}^{(j)} \qquad j \ne i, i+1, \quad
I  \in \mathbb J_\eta \Big |_{i,i+1 \notin I}.
\end{equation}
From (\ref{r1}) and the surrounding text we know that
$s_i c_I(\eta) = c_I(s_i \eta)$ with $I \in \mathbb J_{s_i \eta}$.
We note in general from (\ref{ac1}) that
\begin{equation}\label{bb}
(\overline{s_i \eta})_i =
\left \{ \begin{array}{ll} \bar{\eta}_{i+1}, & \eta_i \ne \eta_{i+1} \\
\bar{\eta}_{i}, & \eta_i = \eta_{i+1} \end{array} \right., \quad
(\overline{s_i \eta})_{i+1} =
\left \{ \begin{array}{ll} \bar{\eta}_{i}, & \eta_i \ne \eta_{i+1} \\
\bar{\eta}_{i+1}, & \eta_i = \eta_{i+1} \end{array} \right., \quad
(\overline{s_{i} \eta})_j = \bar{\eta}_j \quad (j \ne i,i+1).
\end{equation}
Using (\ref{bb}), we see from the definitions (\ref{a2}), (\ref{a3'h})
that for $I \in \mathbb J_\eta \Big |_{i,i+1 \notin I} =
\mathbb J_{s_i \eta} \Big |_{i,i+1 \notin I}$,
$$
A_I(\bar{\eta}/\alpha) = A_I(\overline{s_i \eta}/\alpha), \qquad
\hat{B}_I(\bar{\eta}/\alpha) =
\hat{B}_I(\overline{s_i \eta}/\alpha)
$$
while from (\ref{a3b}) we see that
$$
\tilde{\chi}_I^{(j)}(\bar{\eta}/\alpha) = \tilde{\chi}_I^{(j)}(
\overline{s_i \eta}/\alpha).
$$
Hence (\ref{s5a}) satisfies (\ref{s4y}).

\noindent
{\it The case $i \in I$, $i+1 \notin I$}

We note that in this case $s_i c_I(\eta) \neq c_I(\eta)$.
Consider the first equation in (\ref{s4}) and suppose $\nu = c_I(\eta)$.
We can check from the definitions (\ref{c8}) and (\ref{c8a}) that for
$I \in \mathbb J_\eta \Big |_{i \in I, i+1 \notin I}$,
\begin{equation}\label{mz1}
s_i c_I(\eta) = c_{I \cup \{i+1\}}(\eta)
\end{equation}
so the value of $\tilde{\alpha}_{\eta, s_i \nu}^{(i)}$ on the 
right hand side of
the first equation in (\ref{s4}) is non-zero. Noting from (\ref{u1a})
and the definition (\ref{c8}) of $c_I(\eta)$ that
$$
\bar{\delta}_{i,c_I(\eta)} = (\overline{c_I(\eta)})_i - \bar{\eta}_{i+1},
$$
this equation reads
\begin{equation}\label{mz2}
(1 + \bar{\delta}_{i,\eta}^{-1}) \tilde{\alpha}_{s_i \eta, c_I(\eta)}^{(j)}
= {(\overline{c_I(\eta)})_i - \bar{\eta}_{i+1} -1 \over
(\overline{c_I(\eta)})_i - \bar{\eta}_{i+1}} 
\tilde{\alpha}_{\eta, s_i c_I(\eta)}^{(j)} -
\Big ( {1 \over \bar{\eta}_i - \bar{\eta}_{i+1}} -
{1 \over (\overline{c_I(\eta)})_i - \bar{\eta}_{i+1}} \Big )
\tilde{\alpha}_{\eta, c_I(\eta)}^{(j)}.
\end{equation}

The equation (\ref{mz1}) allows the quantities $A_{I'}$, $\hat{B}_{I'}$
and $\tilde{\chi}^{(j)}_{I'}$, $I' = I \cup \{i+1\}$, making up
$\tilde{\alpha}_{\eta, s_i c_I(\eta)}^{(j)}$ to be related to the
corresponding quantities in $\tilde{\alpha}_{\eta, c_I(\eta)}^{(j)}$.
Thus we can check from the definitions (\ref{a2}) and (\ref{a3'h}) that
\begin{eqnarray}\label{mz2a}
A_{I \cup \{i+1\}}( \bar{\eta}/\alpha) & = &
{(\bar{\eta}_i - (\overline{c_I(\eta)})_i)  \over
(\bar{\eta}_i - \bar{\eta}_{i+1}) ( \bar{\eta}_{i+1} -
 (\overline{c_I(\eta)})_i) } A_I(\bar{\eta}/\alpha) \nonumber \\
\hat{B}_{I \cup \{i+1\}}( \bar{\eta}/\alpha) & = &
{(\bar{\eta}_{i+1} - (\overline{c_I(\eta)})_i) 
\over  \bar{\eta}_{i+1}  - 
(\overline{c_I(\eta)})_i + 1} 
\hat{B}_I( \bar{\eta}/\alpha)\nonumber \\
\tilde{\chi}_{I \cup \{i+1\}}^{(j)} ( \bar{\eta}/\alpha) & = &
\tilde{\chi}_{I}^{(j)} ( \bar{\eta}/\alpha), \qquad j \ne i,i+1.
\end{eqnarray}
When multiplied together according to (\ref{s5a}) to form
$\tilde{\alpha}_{\eta, s_i c_I(\eta)}^{(j)}$ and substituted in
(\ref{mz2}) we find all terms on the right hand side
cancel giving the result
\begin{equation}\label{mz2a'}
\tilde{\alpha}^{(j)}_{s_i \eta, c_I(\eta)} = 0, \qquad j \ne i,i+1.
\end{equation}

Consider now the second and third equation equations in (\ref{s4}) in
the case $\nu = c_I(\eta)$. The requirement in (\ref{s5a}) that
$\tilde{\alpha}^{(j)}_{\eta, c_I(\eta)} \ne 0$ only if $j \in I$, while
we are assuming $i \in I$, $i+1 \notin I$, means the equations read
\begin{eqnarray}\label{mz3}
(1 + \bar{\delta}_{i,\eta}^{-1}) \tilde{\alpha}_{s_i \eta, c_I(\eta)}^{(i+1)}
& = & (1 - \bar{\delta}_{i,\nu}^{-1}) 
\tilde{\alpha}_{\eta, c_{I \cup \{i+1\}}(\eta)}^{(i)} +
\bar{\delta}_{i,\nu}^{-1} \tilde{\alpha}_{\eta, c_{I}(\eta)}^{(i)}
\nonumber \\
(1 + \bar{\delta}_{i,\eta}^{-1}) \tilde{\alpha}_{s_i \eta, c_I(\eta)}^{(i)}
& = & (1 - \bar{\delta}_{i,\nu}^{-1}) 
\tilde{\alpha}_{\eta, c_{I \cup \{i+1\}}(\eta)}^{(i+1)} -
\bar{\delta}_{i,\eta}^{-1} \tilde{\alpha}_{\eta, c_{I}(\eta)}^{(i)} 
\end{eqnarray}
where use has also been made of (\ref{mz1}). To simplify the 
right hand sides of
these equations we note from (\ref{a3b}) that
$$
\tilde{\chi}_{I \cup \{i+1\}}^{(i)}(\bar{\eta}/\alpha) =
{\bar{\eta}_i - \bar{\eta}_{i+1} \over \bar{\eta}_i -
(\overline{c_I(\eta)})_i} \tilde{\chi}_{I}^{(i)}(\bar{\eta}/\alpha),
\quad
\tilde{\chi}_{I \cup \{i+1\}}^{(i+1)}(\bar{\eta}/\alpha) =
{\bar{\eta}_{i+1} - (\overline{c_I(\eta)})_i  \over \bar{\eta}_i -
(\overline{c_I(\eta)})_i} \tilde{\chi}_{I}^{(i)}(\bar{\eta}/\alpha).
$$
Use of these equation, together with the first two equations of (\ref{mz2a})
allows us to express $\tilde{\alpha}^{(i)}_{\eta, c_{I\cup
\{i+1\}}(\eta)}$ and $\tilde{\alpha}^{(i+1)}_{\eta, c_{I\cup
\{i+1\}}(\eta)}$ in terms of $\tilde{\alpha}^{(i)}_{\eta, c_{I}(\eta)}$.
Doing this shows the right hand side
is equal to zero in both cases and so for all
$j=1,\dots,N$
\begin{equation}\label{mz4}
\tilde{\alpha}^{(j)}_{s_i\eta, c_I(\eta)} = 0.
\end{equation}
For $s_i c_I(\eta) \ne c_I(\eta)$
the result (\ref{mz4}) is consistent with (\ref{s5a}) because the result
(\ref{r2}) implies that for $I \in \mathbb 
J_\eta \Big |_{i+1 \notin I, i \in I}$ 
there is no $I' \in    \mathbb J_{s_i \eta}$ such that
$c_{I'}(s_i \eta) = c_I(\eta)$. 

We now proceed to consider the equations (\ref{s4}) in the case
$\nu = s_i c_I(\eta)$, $I \in \mathbb J_\eta
\Big |_{i \in I, i+1 \notin I}$. Proceeding as in the derivation of
(\ref{mz2a'}) we find that in this case the first equation of (\ref{s4})
reads
\begin{eqnarray}\label{su1}
\tilde{\alpha}^{(j)}_{s_i\eta, s_i c_I(\eta)} & = &
{(\bar{\eta}_{i} - \bar{\eta}_{i+1} -1)(\bar{\eta}_{i+1} -
(\overline{c_I(\eta)})_i) \over
(\bar{\eta}_{i} - \bar{\eta}_{i+1})(\bar{\eta}_{i+1} - 
(\overline{c_I(\eta)})_i + 1)} \tilde{\alpha}^{(j)}_{\eta,
 c_I(\eta)} \nonumber \\
& = & 
\tilde{\chi}_{(I \cup \{i+1\})\backslash \{i\}}^{(j)}
(\overline{s_i \eta}/\alpha)
A_{(I \cup \{i+1\})\backslash \{i\}}(\overline{s_i \eta}/\alpha)
\hat{B}_{(I \cup \{i+1\})\backslash \{i\}}(\overline{s_i \eta}/\alpha),
\quad j \ne i,i+1
\end{eqnarray}
where the second equality follows after use of (\ref{s5a}) to
substitute for $\tilde{\alpha}^{(j)}_{\eta,
 c_I(\eta)}$ and use of the definitions (\ref{a2}), (\ref{a3'h}) and
(\ref{a3b}).
An analogous calculation, involving the second and third equations of
(\ref{s4}), gives
\begin{equation}\label{su1a}
\tilde{\alpha}^{(i)}_{s_i\eta, s_i c_I(\eta)} = 0
\end{equation}
as well as the equation (\ref{su1}) in the case $j=i+1$. 
Recalling (\ref{r2})
we see the equations (\ref{su1}) and (\ref{su1a}) are consistent with
(\ref{s5a}).

\noindent
{\it The case $i \notin I$, $i+1 \in I$}

We distinguish the case
\begin{equation}\label{cs1}
s_i c_I(\eta) = c_I(\eta) 
\end{equation}
from 
\begin{equation}\label{cs2}
s_i c_I(\eta) \neq c_I(\eta) 
\end{equation}
In the case (\ref{cs1}) we can check that
\begin{equation}\label{cs2'}
(\overline{ c_I(\eta)})_{i + 1} = \bar{\eta}_i - 1,
\end{equation}
while a feature of the case (\ref{cs2}) is that there is no $I' \in
\mathbb J_\eta$ such that $s_i c_I(\eta) = c_{I'}(\eta)$ and
therefore
\begin{equation}\label{cs3}
\tilde{\alpha}^{(j)}_{\eta, s_i c_I(\eta)} = 0.
\end{equation}

Consider first the equations (\ref{s4}) in the case (\ref{cs1}) (as already
noted, the equations (\ref{s4}) are equivalent to the equations
(\ref{s5}) for $\nu = s_i c_I(\eta) = c_I(\eta)$). The equations read
\begin{eqnarray}\label{ww1}
(1 + \bar{\delta}_{i,\eta}^{-1}) \tilde{\alpha}_{s_i \eta, c_I(\eta)}^{(j)}
&=& (1 - \bar{\delta}_{i,\eta}^{-1}) 
\tilde{\alpha}_{\eta, c_I(\eta)}^{(j)},
\quad j \ne i, i+1 \nonumber \\
(1 + \bar{\delta}_{i,\eta}^{-1}) \tilde{\alpha}_{s_i \eta, c_I(\eta)}^{(i+1)}
&=&  - \bar{\delta}_{i,\eta}^{-1}
\tilde{\alpha}_{\eta, c_I(\eta)}^{(i+1)} \nonumber \\
(1 + \bar{\delta}_{i,\eta}^{-1}) \tilde{\alpha}_{s_i \eta, c_I(\eta)}^{(i)}
& = & \tilde{\alpha}_{ \eta, c_I(\eta)}^{(i+1)}.
\end{eqnarray}
To verify that (\ref{s5a}) satisfies these equations we note that in the
case (\ref{cs1}) the equation (\ref{r5}) is valid, so we should seek
to express (\ref{ww1}) in terms of $\tilde{\alpha}_{s_i \eta,
c_{I \cup \{i\}}(s_i \eta)}$. Now (\ref{bb}), (\ref{a2}), (\ref{a3'h})
and (\ref{a3b}) give
\begin{eqnarray}\label{ww2a}
A_{I \cup \{i\}}(\overline{s_i \eta}/\alpha) & = &
{\bar{\eta}_i - \bar{\eta}_{i+1} - 1 \over
\bar{\eta}_i - \bar{\eta}_{i+1}} A_I(\bar{\eta}/\alpha) \nonumber \\
\hat{B}_{I \cup \{i\}}(\overline{s_i \eta}/\alpha) & = &
{\bar{\eta}_i - \bar{\eta}_{i+1}  \over
\bar{\eta}_i - \bar{\eta}_{i+1} + 1} 
\hat{B}_{I}(\bar{\eta}/\alpha) \nonumber \\
\tilde{\chi}^{(j)}_{I \cup \{i\}}(\overline{s_i \eta}/\alpha) & = & 
\tilde{\chi}^{(j)}_{I}(\bar{\eta}/\alpha), \quad j \ne i, i+1
\nonumber \\
\tilde{\chi}^{(i+1)}_{I \cup \{i\}}(\overline{s_i \eta}/\alpha) & = &
{1 \over \bar{\eta}_{i+1} - \bar{\eta}_i + 1}
\tilde{\chi}^{(i+1)}_{I}(\bar{\eta}/\alpha)  \nonumber \\
\tilde{\chi}^{(i)}_{I \cup \{i\}}(\overline{s_i \eta}/\alpha) & = &
{\bar{\eta}_{i+1} - \bar{\eta}_i \over \bar{\eta}_{i+1} - \bar{\eta}_i + 1}
\tilde{\chi}^{(i+1)}_{I}(\bar{\eta}/\alpha) .
\end{eqnarray}
Making use of these equations in the right hand side of (\ref{ww1}) we find
that for each $j=1,2,\dots,N$ 
\begin{equation}\label{ww3}
\tilde{\alpha}_{s_i \eta, c_I(\eta)}^{(j)} = 
\tilde{\chi}_{I \cup \{i\}}^{(j)}(\overline{s_i \eta}/\alpha)
A_{I \cup \{i\}}(\overline{s_i \eta}/\alpha)
\hat{B}_{I \cup \{i\}}(\overline{s_i \eta}/\alpha)
\end{equation}
which by virtue of (\ref{r5}) is consistent with (\ref{s5a}).

Consider now the equations (\ref{s4}) with $\nu = c_I(\eta)$ in the case
(\ref{cs2}). Then (\ref{c3}) holds, so (\ref{s4}) can be appropriately
simplified. Furthermore, we can check that the second and third members
of (\ref{ww2a}) remain valid, while the remaining equations are to be
replaced by
\begin{eqnarray}\label{aov}
A_{I \cup \{i\}}(\overline{s_i \eta}/\alpha) & = &
- {(\bar{\eta}_{i+1} - (\overline{c_I(\eta)})_{i+1}) \over
(\bar{\eta}_i - \bar{\eta}_{i+1})(\bar{\eta}_i -
(\overline{c_I(\eta)})_{i+1})} A_I(\bar{\eta}/\alpha) \nonumber \\
\tilde{\chi}_{I \cup \{i\}}^{(i+1)}(\overline{s_i \eta}/\alpha) & = &
{\bar{\eta}_i - (\overline{c_I(\eta)})_{i+1} \over
\bar{\eta}_{i+1} - (\overline{c_I(\eta)})_{i+1}}
\tilde{\chi}_{I}^{(i+1)}((\bar{\eta}/\alpha) \nonumber \\
\tilde{\chi}_{I \cup \{i\}}^{(i)}(\overline{s_i \eta}/\alpha) & = &
{\bar{\eta}_{i+1} - \bar{\eta}_i \over
\bar{\eta}_{i+1} - (\overline{c_I(\eta)})_{i+1}}
\tilde{\chi}_{I}^{(i+1)}((\bar{\eta}/\alpha)
\end{eqnarray}
Using these equations to further simplify (\ref{s4}) again 
gives (\ref{ww3}), which we know is consistent with (\ref{s5a}).  
It remains to consider the case $\nu = s_i c_I(\eta)$, 
for which it suffices to restrict attention to the subcase (\ref{cs2}) 
as the subcase (\ref{cs2}) is included in the above working. 
We first simplify the equations (\ref{s4}) according to (\ref{cs3}) 
and then obtain the analogues of (\ref{ww2a}) for the quantities 
$A_{(I \cup \{i\})\backslash\{i+1\}}(\overline{s_i \eta}/\alpha)$
etc.. We find, for $j \ne i+1$, 
$$
\tilde{\alpha}_{s_i \eta, c_I(\eta)}^{(j)} = 
\tilde{\chi}_{I \cup \{i\}\backslash\{i+1\}}^{(j)}(\overline{s_i \eta}/\alpha)
A_{(I \cup \{i\})\backslash\{i+1\}}(\overline{s_i \eta}/\alpha)
\hat{B}_{I \cup \{i\}\backslash\{i+1\}}(\overline{s_i \eta}/\alpha)
$$
while
$$
\tilde{\alpha}_{s_i \eta, c_I(\eta)}^{(i+1)} = 0.
$$
We see from (\ref{r4}) that these equations are consistent with (\ref{s5a}).

\noindent
{\it The case $i,i+1 \in I$}

Analogous to the case $i \notin I$, $i+1 \in I$ we distinguish the case
$s_i c_I(\eta) = c_I(s_i \eta)$ from $s_i c_I(\eta) \notin c_I(s_i \eta)$.
In the latter case
\begin{equation}\label{f2a}
s_i c_I(\eta) = c_{I\backslash \{i+1\}}(\eta).
\end{equation}
This tells us that this case is the same as that with $I' \in 
\mathbb J_\eta$, $i \in I'$, $i+1 \notin I'$, which has already been
dealt with. Thus we can restrict attention to the case
$s_i c_I(\eta) = c_I(s_i \eta)$, when the equations (\ref{s4}) reduce
to the equations (\ref{s5}). Obtaining the analogue of
(\ref{ww2a}) but for $A_{I\backslash \{i\}}(\overline{s_i \eta}/\alpha)$
etc.~we find, for $j \neq i$,
$$
\tilde{\alpha}_{s_i \eta, c_I(\eta)}^{(j)} = 
\tilde{\chi}_{I \backslash\{i\}}^{(j)}(\overline{s_i \eta}/\alpha)
A_{I \backslash\{i\}}(\overline{s_i \eta}/\alpha)
\hat{B}_{I \backslash\{i\}}(\overline{s_i \eta}/\alpha)
$$
while
$$
\tilde{\alpha}_{s_i \eta, c_I(\eta)}^{(i)} = 0.
$$
By virtue of (\ref{r3}) these equations are consistent with
(\ref{s5a}).

\medskip
This completes consideration of the choices of $\nu$ (\ref{cases}) in
all four cases (\ref{Aiv}). In each case it was found (\ref{s5a})
satisfies the recurrences (\ref{s4}), thereby completing the demonstration
that for general $\nu$ (\ref{s5a}) satisfies (\ref{s4}). Since the other
fundamental recurrence (\ref{k5}) has also been shown to be satisfied,
as has the initial condition, our inductive proof is complete. 

\section{An equivalent expansion formula}
\setcounter{equation}{0}
The formula (\ref{fz2}) is the non-symmetric analogue of the Pieri
formula (\ref{pf1}) in the case $p=1$. Here we will use this result
and the formula (\ref{f2}) to derive the analogue of (\ref{pf1}) in the
case $p = N-1$.

First we note that in the case $p=N-1$, analogous to the case $p=1$, the
set $\mathbb J_{N,p}$ appearing in (\ref{f4}) and (\ref{f4'}) can be
indexed by subsets $I = \{t_1,\dots, t_s\}$ of $\{1,\dots,N\}$ with
$t_1 < \cdots < t_s$. Each such subset corresponds to the element
$$
\nu =: \hat{c}_I(\eta)
$$
where
\begin{eqnarray}\label{ptr1}
(\hat{c}_I(\eta))_{t_1} & = & \eta_{t_s}, \nonumber \\
(\hat{c}_I(\eta))_{t_u} & = & \eta_{t_{u-1}} + 1, \quad u=2,\dots,s
\nonumber \\
(\hat{c}_I(\eta))_k & = & \eta_k + 1, \quad k \notin I
\end{eqnarray}
(c.f.~(\ref{c8})). Furthermore, to avoid duplication within the set
$\mathbb J_{N,N-1}$, as with the description (\ref{c8}) of
$\mathbb J_{N,1}$, we must restrict $I$ to maximal subsets 
with respect to $\eta$,
in this case specified by the requirements
\begin{eqnarray}\label{ptr2}
&&\eta_j  \ne  \eta_{t_s} - 1, \qquad j = 1,\dots, t_1 - 1 \nonumber \\
&&\eta_j \ne \eta_{t_u}, \qquad
j=t_u + 1, \dots, t_{u+1} - 1
\end{eqnarray}
for $u=1,\dots,s$ with $t_{s+1} := N + 1$. With this definition of maximal,
analogous to (\ref{c8c}) we define
$$
\hat{\mathbb J}_\eta = \{I: I \: {\rm is \: maximal \: with \: respect
\: to \:} \eta \}.
$$

According to (\ref{f2})
\begin{equation}\label{pq}
c_{\eta \nu}^{(i_1,\dots,i_{N-1})} = c_{\nu \,\eta+(1^N)}^{(j_1)}
{\ml E_\eta | E_\eta \mg \over \ml E_\nu | E_\nu \mg}.
\end{equation}
Now $c_{\nu \, \eta + (1^N)}^{(j_1)}$ is non-zero only if
$\eta + (1^N) = c_I(\nu)$, $I \in \mathbb J_{\nu}$,
$j_1 \in I$. With $I = \{t_1,\dots, t_s\}$ we see from (\ref{c8}) that
$\eta + (1^N) = c_I(\nu)$ gives
\begin{eqnarray}
\eta_{t_j} & = & \nu_{t_{j+1}} - 1, \qquad j=1,\dots, s-1 \nonumber \\
\eta_{t_s} & = & \nu_{t_1} \nonumber \\
\eta_i & = & \nu_i - 1, \qquad i \notin I
\end{eqnarray}
while from (\ref{c9}) the condition $I \in \mathbb J_{\nu}$
gives
\begin{eqnarray}
\eta_j & \ne & \eta_{t_s} - 1 \qquad j=1,\dots, t_1 - 1 \\
\eta_j & \ne & \eta_{t_u} \qquad j=t_u+1, \dots, t_{u+1} - 1
\end{eqnarray}
for $u=1,\dots,s$ with $t_{N+1} := N+1$. These are precisely the equations
(\ref{ptr1}) with $\nu := \hat{c}_I(\eta)$ and the equations (\ref{ptr2})
for $I \in  \hat{\mathbb J}_{\eta}$, so we conclude
$$
\eta = 
c_I(\nu) \Big |_{I \in \mathbb J_{\nu}} \quad {\rm iff} \quad 
\nu + (1)\hat{c}_I(\eta) \Big |_{I \in \hat{{\mathbb  J}}_\eta}.
$$

It remains to substitute for the explicit values in the right hand side
of (\ref{pq}).
With $\eta + (1^N) = c_I(\nu)$, $I \in \mathbb J_{\nu}$
and thus $\nu = \hat{c}_I(\eta)$, $I \in 
\hat{\mathbb J}_\eta$  we read off from (\ref{obt}) and
(\ref{en}) that
\begin{equation}\label{fii}
c_{\eta \nu}^{(i_1,\dots,i_{N-1})}  = 
{e_\eta d_\nu \over d_\eta e_\nu}
{d_\eta e_{\eta + (1^N)} \over d_{\eta + (1^N)} e_\eta}
\tilde{\chi}^{(j_1)}_I(\overline{\eta}/\alpha)
A_I(\overline{\eta}/\alpha)
\hat{B}_I(\overline{\eta}/\alpha)
\end{equation}
where use has been made of the facts that ${\cal N}_\eta =
{\cal N}_{\eta + (1^N)}$ and $\tilde{\chi}_I^{(j_1)}
((\overline{\eta} + c)/\alpha)
= \tilde{\chi}^{(j_1)}_I(\overline{\eta}/\alpha)$ etc.~for any constant
$c$.
This further simplifies by noting from (\ref{d1}), (\ref{leg}) and (\ref{ac1})
that
$$
{d_\eta \over d_{\eta + (1^N)}}
= {1 \over \prod_{j=1}^N(\bar{\eta}_j + \alpha + N)}
$$
while (\ref{e1}) together with (\ref{bb}) and (\ref{ac1}) implies
$$
{e_{\eta + (1^N)} \over e_\eta} =
\prod_{j=1}^N(\bar{\eta}_j + \alpha + N).
$$
Substituting these formulas in (\ref{fii}) gives
\begin{equation}\label{56}
c_{\eta \nu}^{(i_1,\dots,i_{N-1})}  =
{e_\eta d_\nu \over d_\eta e_\nu}
\tilde{\chi}^{(j_1)}_I(\overline{\eta}/\alpha)
A_I(\overline{\eta}/\alpha)
\hat{B}_I(\overline{\eta}/\alpha), \quad
\nu = \hat{c}_I(\eta), \: I \in \hat{\mathbb J}_\eta
\end{equation}
(c.f.~(\ref{s5b'})).

As in the derivation of (\ref{fu2}) from (\ref{fu1}), if follows from
(\ref{56}) that
\begin{equation}\label{57}
e_{N-1}(z) E_\eta(z) = - \alpha {e_\eta \over d_\eta}
\sum_{I \in \hat{\mathbb J}_\eta}
{e_{\hat{c}_I(\eta)} A_I(\bar{\eta}/\alpha) \hat{B}_I(\bar{\eta}/\alpha)
\over d_{\hat{c}_I(\eta)} } E_{\hat{c}_I(\eta)}(z).
\end{equation}

\section{The coefficient $A_{\eta, \nu}^{(p)}$ of the Pieri type formula
for general $p$}
\setcounter{equation}{0}
In this final section we will consider features of the coefficient
$A_{\eta, \nu}^{(p)}$ in the expansion (\ref{f4'}) for general $p$.
Our first result concerns the value of $A_{\eta, \nu}^{(p)}$ for a
particular value of $\nu$. Now, according to (\ref{2.8'}) the coefficient
of $z^\eta$ in $E_\eta(z)$ is unity and all other monomials are smaller with
respect to the ordering $\triangleleft$. Consider all sets of the form
$M = \{j_1,\dots,j_p | 1 \le j_1 < \cdots < j_p \le N \}$ and let
$$
(\chi_M)_i = \left \{ \begin{array}{ll}1, & i \in M \\
0, & {\rm otherwise} \end{array} \right.
$$
so that the $p$th monomial symmetric function can be written
$$
e_p(z) = \sum_M z^{\chi_M}.
$$
Let $M^*$ be the particular set $M$ such that
$$
\eta + \chi_M \triangleleft \eta + \chi_{M^*}
$$
for all $M \ne M^*$. Then we must have
\begin{equation}\label{aw0}
A_{\eta, \eta + \chi_{M^*}}^{(p)} = 1.
\end{equation}
Moreover, with
$$
l_\eta'(i) := \#\{k < i | \eta_k \ge \eta_i \} + \#\{k > i| \eta_k > \eta_i
\}
$$
it follows from the definition of $\triangleleft$ that
\begin{equation}\label{aw}
( \eta + \chi_{M^*})_i = \left \{
\begin{array}{ll} \eta_i + 1, & l_\eta'(i) \le p-1 \\ 
\eta_i, & l_\eta'(i) \ge p \end{array} \right.
\end{equation}

Associated with (\ref{aw}) are the sets
\begin{equation}
G_0 := \{i \in \{1,\dots,N\}: l_\eta'(i) \ge p \}, \qquad
G_1 := \{i \in \{1,\dots,N\}: l_\eta'(i) \le p - 1 \}.
\end{equation}
An alternative characterization follows by noting that since $\eta
\subseteq \eta + \chi_{M^*}$ we have $\eta \preceq \eta + \chi_{M^*}$
and so from the definition of $\preceq$, for $\nu =
 \eta + \chi_{M^*}$,
\begin{equation}
G_0 := \{i \in \{1,\dots,N\}: \nu_{\pi(j)} = \eta_j \}, \quad
G_1 := \{i \in \{1,\dots,N\}: \nu_{\pi(j)} = \eta_j + 1 \}.
\end{equation}
Let us now put 
\begin{equation}\label{siv}
B_{\eta,\nu}^{(p)} := {d_\nu' e_\eta' \over e_\nu' d_\eta'}
A_{\eta, \nu}^{(p)}.
\end{equation}
Then it follows from the definitions (\ref{d1}) and (\ref{e1}), together
with (\ref{aw0}), that
\begin{equation}\label{swa}
B_{\eta, \eta + \chi_{M^*}}^{(p)} =
\Big ( \prod_{j \in G_0, \, k \in G_1 \atop j<k}
{\bar{\eta}_j - \bar{\eta}_k +1 \over \bar{\eta}_j - \bar{\eta}_k}
\Big )
\Big ( \prod_{j \in G_1, \, k \in G_0 \atop \pi(j)<\pi(k)}
{\bar{\eta}_j - \bar{\eta}_k +\alpha - 1 \over 
\bar{\eta}_j - \bar{\eta}_k + \alpha}
\Big )
\end{equation}
In the case $p=1$, comparison with the expression
\begin{equation}\label{san}
B_{\eta,\nu}^{(1)} = - \alpha A_I(\bar{\eta}/\alpha) \hat{B}_I
(\bar{\eta}/\alpha), \quad I \in {\mathbb J}_\eta,
\end{equation}
which follows from (\ref{obt}), (\ref{4.47'}) and (\ref{siv}), we see
that as written (\ref{swa}) is in fact valid for all $\nu = c_I(\eta)$
with $I$ consisting of a single element $t_1$. More explicitly, we
then have
\begin{eqnarray}
&&
A_I(\bar{\eta}/\alpha) = a(\bar{\eta}_{t_1}/\alpha - 1, 
\bar{\eta}_{t_1}) = - 
{1 \over \alpha} \label{san1}  \\
&&\hat{B}_1(\bar{\eta}/\alpha) =
\prod_{j=1}^{t_1 - 1} b(\bar{\eta}_{t_1}/\alpha,
\bar{\eta}_j/\alpha)
\prod_{j=t_1 + 1}^N b(\bar{\eta}_{t_1}/\alpha+1,
\bar{\eta}_j/\alpha);  \label{san2}
\end{eqnarray}
the factor (\ref{san1}) cancels with $-\alpha$ in (\ref{san}) while the
two products (\ref{san2}) correspond with the two products in (\ref{swa})
respectively.

We can extend the form (\ref{swa}) so that for $p=1$ there is agreement
with (\ref{san}). The extended form is
\begin{eqnarray}
B_{\eta,\nu}^{(p)} & = & \Big (
\prod_{j \in G_0, k \in G_1 \atop j < k}
{\bar{\eta}_j - \bar{\eta}_k + 1 \over
\bar{\eta}_j - \bar{\eta}_k } \Big )
\Big ( \prod_{j \in G_1, \, k \in G_0 \atop \pi(j)<\pi(k)}
{\bar{\eta}_j - \bar{\eta}_k +\alpha - 1 \over
\bar{\eta}_j - \bar{\eta}_k + \alpha}
\Big ) \nonumber \\
&& \times
\prod_{\pi^2(j) < \pi(j) < j}
{1 \over \bar{\eta}_{\pi(j)} - \bar{\eta}_j}
\prod_{j \le \pi^2(j) \le \pi(j)}
{1 \over \bar{\eta}_{\pi(j)} - \bar{\eta}_j - \alpha},
\end{eqnarray}
valid for $p=1$ and $\nu = c_I(\eta)$, $I \in \mathbb J_\eta$.
The significant feature of this form is that in the general $p$ case, with
$\nu = c_I(\eta) \in \mathbb J_{N,p}$ and $I$ such that at most one part
of $\eta$ in the formation of $\nu$ according to the prescription
below (\ref{c7}) move downwards, explicit small $N$ calculations
indicate it remains valid. However we have no proof of this empirical
observation.

\section*{Acknowledgements}
We thank Dan Marshall for providing us with his notes on the derivation
of (\ref{fu1}). This work was supported by the Australian Research
Council.


\end{document}